\documentclass[12pt]{article}

\usepackage{amsmath,amsfonts,amssymb,amsthm,newlfont}
\newcommand{\bb}{\begin{equation}}
\newcommand{\ds}{\displaystyle }
\newcommand{\ee}{\end{equation}}

\def\ga{\gamma}

\def\al{\alpha}

\def\be{\beta}

\def\l1{{\lambda}_1}

\def\nb{\nabla}
\def\kd{\partial}

\def\nb{\nabla}

\begin{document}
\pagenumbering{arabic}

\title{\huge \bf A Caffarelli-Kohn-Nirenberg type inequality
on Riemannian Manifolds}
\author{\rm \large Yuri Bozhkov \\ \\
\it Instituto de Matem\'atica, Estatistica e \\
\it Computa\c c\~ao Cient\'\i fica - IMECC \\
\it Universidade Estadual de Campinas - UNICAMP \\
\it C.P. $6065$, $13083$-$970$ - Campinas - SP, Brasil \\
\rm E-mail: bozhkov@ime.unicamp.br }
\date{15th June 2009}
\maketitle

\begin{abstract}
\rm We establish a generalization to Riemannian manifolds of the
Caffarelli-Kohn-Nirenberg inequality. The applied method is based
on the use of conformal Killing vector fields and Enzo Mitidieri's
approach to Hardy inequalities.
\end{abstract}

%\vskip 1cm
\begin{center}
{2000 AMS Mathematics Classification numbers: 58E35, 26D10}
\end{center}

%\vskip 1cm
\begin{center}
\end{center}

\pagenumbering{arabic}
\parindent 30pt
\baselineskip 16pt
%\newpage

\section{Introduction}

\

It is well-known that Hardy type inequalities have been widely
used in the theory of differential equations. An inequality of
this kind is the celebrated general interpolation inequality with
weights obtained by Caffarelli, Kohn and Nirenberg in \cite{ckn}.
The latter can be stated in the form of the following

\

{\bf Theorem 1.} \rm (Caffarelli-Kohn-Nirenberg \cite{ckn}) {\it
For all $a,b\in \mathbb{R}$ and $u\in C^{\infty }_0
(\mathbb{R}^n\setminus \{0\})$ the inequality
 \bb\label{i1} C \int_{\mathbb{R}^n} \frac{|u|^2}{|x|^{a+b+1}}\;dx \leq
 \left ( \int_{\mathbb{R}^n} \frac{|u|^2}{|x|^{2a}}\;dx \right )^{1/2}
 \left ( \int_{\mathbb{R}^n} \frac{|\nb u|^2}{|x|^{2b}}\;dx \right )^{1/2}\ee
 holds with sharp constant $C=C(a,b)= |n-(a+b+1)|/2$.}

 \

 \rm (We note that the presented
 above Theorem 1 is the version of the Caffarelli-Kohn-Nirenberg
 Inequality discussed in \cite{dc}.)

 There is a very big number of recent papers dedicated to
 generalizations of the Hardy, Sobolev and Caffarelli-Kohn-Nirenberg inequalities in various
 contexts.
 For this reason we merely cite just a few of them (e.g.
 \cite{acp,agg,cx,car,cw,dc,em2,x})
 directing the interested reader to use an internet search
 instrument for further references.

 Since most of the problems in differential geometry can be reduced to
 problems in differential equations on Riemannian manifolds, it is important
 to have in this case investigation tools similar to those successfully used
 in the Euclidean case, in particular, integral inequalities for functions in
 weighted Sobolev type spaces. Hence the purpose of this work is to obtain another
 generalization of (\ref{i1}) to Riemannian manifolds.

 Let $M$ be a complete $n-$dimensional manifold, $n\geq
3$, endowed with a Riemannian metric $g=(g_{ij})$ given in local
coordinates $x^i$, $i=1,2,...,n$, around some point. We suppose
that on $M$ exist (nontrivial) conformal Killing vector fields
$h=h^i\frac{\kd}{\kd x^i}$:
 \bb\label{a0} {\nb }^i h^j + {\nb }^j h^i= \frac{2}{n}
 (div\; h)\;g^{ij} =\mu\; g^{ij}. \ee
Here ${\nabla }^i$ is the covariant derivative corresponding to
the Levi-Civita connection, uniquely determined by $g$, $div\; h$
is the covariant divergence operator and summation
 from $1$ to $n$ over repeated Latin indices is understood.

 Now we state the main result of this paper:

 \

 {\bf Theorem 2.} {\it Let $p>1$ and $u\in W^{1,p}(M)$. Suppose that $M$ admits
 a $C^1$ conformal Killing vector field $h=h^i\frac{\kd}{\kd x^i}$ such
 that $div\;h >0$. Then}
 \bb\label{har1} C \int_M
 \frac{div\;h}{|h|^{a+b+1}}\;|u|^p\;dV\leq
 \left (\int_M \frac{div\;h}{|h|^{aq}}\;|u|^p \;dV \right )^{1/q}
 \left (\int_M \frac{(div\;h)^{1-p}}{|h|^{bp}}\;|\nb u|^p \;dV \right )^{1/p} \ee
 {\it where $dV$ is the volume form of $M$,} $|h|^2 =g_{ij}h^ih^j$, $\frac{1}{p} +\frac{1}{q}=1$ {\it and the
 constant $C=|n-(a+b+1)|/(p n)$ is sharp.}

 \

 The proof of this result is based on the idea of the E. Mitidieri's vector field approach to Hardy
 inequalities, proposed in
 \cite{em2}, which is applied here to conformal Killing vector fields on Riemannian
 manifolds.

 Note that in Theorem 2 we do not impose additional
 conditions neither on the gradient of the weight function nor on
 its second derivatives (e.g. certain kind of super-harmonicity).

 Further we assume that $M$ admits a conformal Killing vector field
$h = {h }^i \frac{\kd }{\kd x^i}$ such that \bb\label{u2} {\nabla
}^i {h }^j + {\nabla }^j {h}^i = c \;g^{ij} =\frac{2}{n}\; (div\;h
) \;g^{ij},  \ee
 where $c\neq 0$ is a constant. This means
we suppose that $M$ admits a homothety (`homothetic motion') which
is not an infinitesimal isometry of $M$. Actually we may choose
$c=2$ in (\ref{u2}) and hence $div\; h = n$ without loss of
generality. (Otherwise, since $c\neq 0$, we could consider $2 h
/c$ instead of $h $.) Such a special conformal vector field will
be used to obtain the following {\it exact} generalization of the
 Caffarelli-Kohn-Nirenberg Inequality in ${\mathbb{R}}^n$.

 \

 {\bf Corollary 1.} {\it If $div\;h=n$ and $p=q=2$, then
 $(\ref{har1})$ becomes}
 \bb\label{i3} C \int_{M} \frac{|u|^2}{|h|^{a+b+1}}\;dV \leq
 \left ( \int_{M} \frac{|u|^2}{|h|^{2a}}\;dV \right )^{1/2}
 \left ( \int_{M} \frac{|\nb u|^2}{|h|^{2b}}\;dV \right )^{1/2}\ee
 {\it with sharp constant $C=C(a,b)= |n-(a+b+1)|/2$.}

 \

 Note that for $M={\mathbb{R}}^n$, $g^{ij}={\delta }^{ij}-$the Euclidean
 metric of ${\mathbb{R}}^n$ and $h =x^i \frac{\kd }{\kd x^i}$, the radial vector field
 corresponding to a dilational transformation in ${\mathbb{R}}^n$, the inequality
 (\ref{i3}) reduces to the Caffarelli-Kohn-Nirenberg Inequality (\ref{i1}).

 If $a=p-1$, $b=0$ in (\ref{har1}), then we can get easily the sharp
 Hardy type
 inequality established in \cite{ye1}.
 Namely:

 \

 {\bf Corollary 2.} {\rm (\cite{ye1})} {\it Let $p>1$ and $u\in W^{1,p}(M)$. Suppose that $M$
 admits a $C^1$
 conformal Killing vector field $h$ such that $div\;h >0$. Then
 $(div\;h)|u|^p/|h|^p\in L^1(M)$
 and }
 \[ \left (\frac{|n-p|}{n p}\right )^p \int_M
 \frac{div\;h}{|h|^p}\;|u|^p\;dV\leq \int_M (div\;h)^{1-p}\;|\nb u|^p \;dV.
 \]

 \

 We observe that the Caffarelli-Kohn-Nirenberg Inequality (\ref{i1}) with
 $b=0$ and $a=-1$ represents the classical Uncertainity Principle
 of quantum mechanics:
 \[ \frac{n^2}{4}\left ( \int_{\mathbb{R}^n} |u|^2dx \right )^2 \leq
 \left ( \int_{\mathbb{R}^n} |x|^2 |u|^2 dx \right )
 \left ( \int_{\mathbb{R}^n} |\nb u|^2dx \right ).\]
 Then Theorem 2 immediately implies its Riemannian analog:

 \

 {\bf Corollary 3.} {\it Let $p>1$, $1/p+1/q=1$ and
 $u\in C^{\infty }_0 (M\setminus \{zeros\;\; of\;\;h\})$.
 Suppose that $M$ admits a $C^1$ conformal Killing vector field $h$
 such that $div\;h >0$. Then the following Uncertainity Principle
 inequality holds:}
 \[\frac{1}{p} \int_M (div\;h) \;|u|^p\;dV\leq
 \left (\int_M (div\;h) \;|h|^{q}\;|u|^p \;dV \right )^{1/q}
 \left (\int_M (div\;h)^{1-p}\;|\nb u|^p \;dV \right )^{1/p}.  \]

 \

 Another result which can be handled by the technique developed in this paper
 is one of the results obtained
 in \cite{x}. In that work C. Xia proves, among other things, the following sharp
 Caffarelli-Kohn-Nirenberg inequality for any $u\in C^{\infty }_0
(\mathbb{R}^n)$:
 \[ \int_{\Bbb R^n}|x|^{\gamma r}|
 u|^r\,dx \leq {r\over n+\gamma r}\left(\int_{\Bbb R^n}|
 x|^{\alpha p}|\nabla u|^p\,dx\right)^{1/p}
 \left(\int_{\Bbb R^n}| x|^\beta|
 u|^{p(r-1)\over p-1}\,
 dx\right)^{(p-1)/p}, \]
 where $n\ge2$, $1<p<r$, $\alpha>0$, $\beta>0$ satisfy
 \bb\label{cx} {1\over p}+{\alpha\over n}>0,\;\;\; {p-1\over p(r-1)}+{\beta\over
 n}>0,\;\;\;
 {1\over r}+{\gamma\over n}>0\;\;\;\gamma={\alpha-1\over r}+{(p-1)\beta\over pr}.
 \ee

 The method used in the present paper permits to obtain its
 generalization:

 \

 {\bf Corollary 4.} {\it Suppose that the parameters $\al ,\be , \gamma , r, p$
 satisfy the condition $(\ref{cx})$ and the function $u\in C^{\infty }_0(M)$.
 In addition suppose that $M$ admits a $C^1$ conformal Killing vector field $h$
 such that $div\;h >0$. Then}
 \[\int_M (div\;h)|h|^{\ga r} \;|u|^r\;dV\leq \]\bb\label{az4} \ee\[
 \leq\frac{r n}{n+\ga r} \left (\int_M
 (div\;h)^{1-p}\;|h|^{\al p}\;|\nb u|^p \;dV \right )^{1/p}
 \left (\int_M (div\;h) \;|h|^{\be }\;|u|^{(r-1)q} \;dV \right
)^{1/q},  \]
 {\it where } $q= p/(p-1).$

 \

 We shall prove Theorem 2 and Corollary 4 in the next section. In
 section 3 we briefly present comments and concluding
 remarks.

 \

  \section{An application of the vector fields method}

 \

 We begin with a preliminary

 \

 {\bf Lemma.} {\it If $h$ is a conformal Killing vector
 field satisfying $(\ref{a0})$ then }
 \bb\label{har2} div \left ( \frac{h}{|h|^{k}}\right )=
 \frac{n-k}{2}\frac{\mu }{|h|^{k} },\ee
 {\it where} $k\in\mathbb{R} $. {\it If } $k>0$, {\it the relation
 $(\ref{har2})$ holds on $M\setminus \{ zeros\;\; of\;\; h \}$;
 otherwise, it holds on the whole of $M$.}

 \

 {\bf Proof.} From (\ref{a0}) we have that the covariant
 divergence
 \bb\label{a1} div \left ( \frac{h}{|h|^{k}}\right )=
 \frac{n}{2}\frac{\mu}{|h|^{k}}-k|h|^{-k-2}h_jh^i {\nb }_i
 h^j.\ee
 Further we interchange the indices $i$ and $j$ and use
 (\ref{a0}) in the following way:
 \[h_j h^i {\nb }_i h^j = h_j h_i {\nb }^i h^j =h_j h_i (-{\nb }^j h^i +\mu\; g^{ij})=
 - h_j h^i {\nb }_i h^j +\mu\;|h|^2 .\]
 Hence
 \bb\label{a2} h_j h^k {\nb }_k h^j =\frac{\mu }{2} |h|^2. \ee
 Thus the relation (\ref{har2}) follows from (\ref{a1}) and
 (\ref{a2}).

 \

 After this little preparatory work we are now ready for the

 \

 {\bf Proof of Theorem 4.} As usual, we may
 assume that $u\in C^{\infty }_0(M)$ without loss of generality. Then integrating by parts
 and using (\ref{har2}) with $k=a+b+1$ we obtain that
 \[\ds{\int_M}
\ds{\frac{|u|^{p-1}(h,\nb u) }{|h|^{a+b+1}}}\;dV=
-\frac{n-(a+b+1)}{2 p} \ds{\int_M} \ds{\frac{\mu
}{|h|^{a+b+1}}}|u|^p\;dV .\] Hence
\[ \frac{|n-(a+b+1)|}{2 p} \ds{\int_M} \ds{\frac{\mu }{|h|^{a+b+1}}}|u|^p\;dV
\leq \ds{\int_M} \ds{\frac{|u|^{p-1}|\nb u |}{|h|^{a+b}}}\;dV
\]

\[= \ds{\int_M} \ds{\frac{|u|^{p/q} {\mu }^{1/q}}{|h|^{a}}}.
\frac{|\nb u |{\mu }^{-1/q}}{|h|^{b}}\;dV
\]

\[\leq \left ( \ds{\int_M} \ds{\frac{\mu }{|h|^{aq}}}|u|^p \;dV\right )^{1/q}\;
\left ( \ds{\int_M} \ds{\frac{{\mu }^{1-p} }{|h|^{bp}}}|\nb u|^p
\;dV\right )^{1/p} \] by the H\"older inequality with $q=p/(p-1)$.
Thus \[  \frac{|n-(a+b+1)|}{2 p} \ds{\int_M} \ds{\frac{\mu
}{|h|^{a+b+1}} |u|^p }\;dV\leq \left ( \ds{\int_M} \ds{\frac{\mu
}{|h|^{aq}}}|u|^p \;dV\right )^{1/q}\; \left ( \ds{\int_M}
\ds{\frac{{\mu }^{1-p} }{|h|^{bp}}}|\nb u|^p \;dV\right )^{1/p}.\]
Then (\ref{har1}) follows from the latter inequality with $\mu =2
(div\;h)/n $ (see (\ref{a0})). It is clear that the constant $C$
is sharp \cite{dc,em2}.

\

{\bf Proof of Corollary 4.} \rm By (\ref{har2}) with $k=-\ga r$ we
have
\[\ds{\int_M}
\ds{|h|^{\ga r}|u|^{r-1}(h,\nb u)}\;dV= -\frac{n+\ga r}{2 r}
\ds{\int_M} \ds{\mu |h|^{\ga r}|u|^r}\;dV \] after integration by
parts. Then
\[ \ds{\int_M} \ds{\mu |h|^{\ga r}}|u|^r\;dV
\leq \frac{2 r}{n+\ga r}\ds{\int_M} \ds{|h|^{\ga r +
1}|u|^{r-1}|\nb u |}\;dV
\]

\[= \frac{2 r}{n+\ga r}\ds{\int_M} \ds{ |h|^{\al }|\nb u |{\mu }^{-1/q}}.
{\mu }^{1/q}|h|^{\ga r +1-\al}|u|^{r-1}\;dV \]

\[\leq  \frac{2 r}{n+\ga r}\left (\ds{\int_M} {\mu }^{1-p}\;|h|^{\al p}\;|\nb u|^p \;dV \right )^{1/p}
 \left (\int_M {\mu } \;|h|^{(\ga r +1-\al )q }\;|u|^{(r-1)q} \;dV \right
)^{1/q} \] by the H\"older inequality with $q=p/(p-1)$. Hence we
obtain (\ref{az4}) setting $\be = (\ga r +1-\al)p/(p-1)$ which is
equivalent to the last relation in (\ref{cx}).

 \

 \section{Conclusion}

 \

 We have proved a Riemannian analog of the Caffarelli-Kohn-Nirenberg
 Inequality using the Enzo Mitidieri's approach to Hardy inequalities, devised and
 developed in \cite{em2}, applied to conformal Killing vector
 fields. We would also like to observe that the presented proof of
 Theorem 2 can be considered as another proof of the Caffarelli-Kohn-Nirenberg
 Inequality in ${\mathbb{R}}^n$ simply taking $h =x^i \frac{\kd }{\kd x^i}$,
 the radial vector field which has been used in \cite{dc,em2}.

 We point out that the quadratic inequality
 method used by David Costa in his elegant paper \cite{dc} applies to
 Riemannian manifolds as well. Indeed, if $h= {h }^i \frac{\kd }{\kd x^i}$ is a conformal
 Killing vector field such that $div\;h=n$, consider another vector field
 \[ W = \frac{\nb u}{|h|^b} +
 t\;u\frac{h}{|h|^{a+1}}=\left ( \frac{1}{|h|^b} g^{ij} u\frac{\kd u}{\kd x^i} +
 t\;u\frac{h^i}{|h|^{a+1}} \right )\frac{\kd }{\kd x^i}
 \]
 where $t\in\mathbb{R}$ and $g^{ij} $ is the inverse matrix of $g_{ij} $. From $g(W,W)\geq 0$,
 that is, from $ g_{ij}W^i W^j\geq 0$, and the identity (\ref{har2}), it is
 obvious
 that one can prove (\ref{i3}) following the argument in
 \cite{dc}, p. 313.

 \

 \section*{Acknowledgements}

 \

 \rm We would like to thank Professors Enzo Mitidieri and David
 Costa for their encouragement. We would also like to thank CNPq,
 Brasil, for partial
 financial support.

 \end{document}